\documentclass[12pt,a4paper,draft]{article}
\usepackage[T2A]{fontenc}
\usepackage[english]{babel}
\usepackage{latexsym,amsfonts,amssymb,amsmath,longtable,amsthm,mathrsfs,cite}
\renewcommand{\le}{\leqslant}
\newcommand{\prn}{\,\,{\rm prn}\,\,}

\renewcommand{\qed}{\hfill{$\Box$}}

\newtheorem*{Prop*}{{\bfseries Proposition}}
\newtheorem*{Cor*}{{\bfseries Corollary}}

\newtheorem{Lemma}{{\bfseries Lemma}}
\newtheorem{Cor}[Lemma]{{\bfseries Corollary}}
\newtheorem{Theo}[Lemma]{{\bfseries Theorem}}
\theoremstyle{definition}

\DeclareMathOperator{\GL}{GL}

\title{\vspace{-1cm} \hfill{\normalsize MSC2010 20D20}{
\fontfamily{cmr} \fontseries{bx} \selectfont \\ \vspace{1cm} The existence of pronormal $\pi$-Hall subgroups in  $E_\pi$-groups}}
\author{\bf  D. O. Revin, E. P. Vdovin}

\begin{document}
\sloppy


\maketitle

To Yu.L.Ershov on his seventy-fifth birthday

\begin{abstract}

A subgroup $H$ of a group $G$ is called {\it pronormal},
if the subgroups  $H$ and $H^g$ are conjugate in  $\langle H, H^g\rangle$  for 
every  $g\in G$. It is proven that if a finite group
$G$ possesses a $\pi$-Hall subgroup for a set of primes
$\pi$, then its every normal subgroup (in particular,  $G$ itself) has a
$\pi$-Hall subgroup that is pronormal in~$G$.
\end{abstract}


\section*{Introduction}

Throughout the paper  the term ``group'' implies
``finite group''.

According to the definition by P.~Hall, a subgroup $H$ of a group $G$ is called {\it
pronormal}, if the subgroups  $H$ and $H^g$ are conjugate in~$\langle H, 
H^g\rangle$  for every~${g\in G}$.

The pronormality for subgroups is a more general property, than normality, and 
it plays
an important role in the theory of groups. In particular,
the Frattini argument holds for pronormal subgroups: {\sl If a pronormal 
subgroup  $H$ of $G$ lies in a
normal subgroup~$A$, then $G=AN_G(H).$} This statement is often  important
in inductive arguments. Notice that    $G=AN_G(H)$ if and only if
$ 
H^G=\{H^g\mid g\in G\}$
 and
$H^A=\{H^a\mid a\in A\}$ coincide.

In view of the  Sylow theorem, the  Sylow subgroups of a group, sa well as the  
Sylow subgroups
of every normal subgroup are examples of pronormal subgroups. The goal of the 
paper is
to study in what form this properties of Sylow subgroups can be reformulated for
Hall subgroups. We recall the   appropriate definitions.

Throughout the paper we suppose that $\pi$ is a fixed set of primes. We denote
by  $\pi'$
the set of all primes not in $\pi$; and by
$\pi(n)$, the set of all prime divisors of a natural number  $n$, while
for a group  $G$ we denote the set   $\pi(|G|)$  by  $\pi(G)$.
A natural   $n$ with $\pi(n)\subseteq\pi$ is called a
{\it $\pi$-number}, while a group $G$ with $\pi(G)\subseteq \pi$ is called a
{\it $\pi$-group.}
A subgroup $H$ of $G$ is called a {\it $\pi$-Hall} subgroup, if
$\pi(H)\subseteq\pi$ and $\pi(|G:H|)\subseteq \pi'$. Thus in case
$\pi=\{p\}$ the definition of  $\pi$-Hall subgroup coincides with the usual notion
of Sylow  $p$-subgroup. A subgroup is said to be a {\em Hall} subgroup, if it 
is a
$\pi$-Hall subgroup for a set of primes~$\pi$, i.~e. if its index and  order are
coprime.

According to \cite{Hall} we say that  $G$ {\it satisfies $E_\pi$} (or briefly
$G\in E_\pi$), if $G$ possesses a $\pi$-Hall subgroup. If, at that, every two 
$\pi$-Hall
subgroups are conjugate, then we say that $G$ {\it satisfies
$C_\pi$} ($G\in C_\pi$). The group satisfying
$E_\pi$ or $C_\pi$ we call also an $E_\pi$- or a $C_\pi$-{\it group} respectively.

The Hall theorem implies that Hall subgroups are pronormal in solvable groups.
Also the $\pi$-Hall subgroups are known to be pronormal
\begin{itemize}
 \item in finite simple groups~\cite{VR3};
\item  in $C_\pi$-groups~\cite{VR4}.
\end{itemize}
In~\cite{GR} some properties were specified of the groups such that $\pi$-Hall 
subgroups exist and all are pronormal.

At the same time it is known~\cite{VR4}
that if a set  $\pi$ of primes is such that
$E_\pi\ne C_\pi$, then for every  $X\in E_\pi\setminus C_\pi$ and  $p\in\pi'$
the group $X\wr \mathbb{Z}_p$ possesses nonpronormal
$\pi$-Hall subgroups, so we cannot transfer the property of pronormality
of Sylow subgroups to Hall subgroups.

The following analog of the Frattini argument for $\pi$-Hall subgroups is
proven in~\cite[Theorem~1]{VR5}: {\sl If $G\in E_\pi$, then each normal 
subgroup $A$
of $G$ possesses a $\pi$-Hall subgroup  $H$ such that  $G=AN_G(H)$.}

In the paper we prove the following statement, on using results from
\cite{RevVdoArxive,ExCrit,VR1,VR2,VR3,VR5}.

\begin{Theo}\label{MainTheorem} Let $G\in E_\pi$ for a set of primes
$\pi$ and $A\trianglelefteq  G$. Then there is a $\pi$-Hall subgroup of $A$
pronormal in $G$.
\end{Theo}

\begin{Cor}\label{cor}
Let $\pi$ be a set of primes. Then every   group possessing a
$\pi$-Hall subgroup has a pronormal  $\pi$-Hall subgroup.
\end{Cor}

Notice that the result of~\cite{VR2} on the pronormality of $\pi$-Hall subgroups in
$C_\pi$-groups (or, equivalently, on the inheriting of the $C_\pi$-property by
overgroups of $\pi$-Hall subgroups) is a particular case of this statement.

Theorem~\ref{MainTheorem} generalizes the above mentioned result of~\cite{VR5}.
Notice also the following statement that generalizes the  useful 
Lemma~\ref{HallExist} (located below)
and giving a criterion of existence of  $\pi$-Hall subgroups in nonsimple groups.

\begin{Cor}\label{cor1} Let $A\trianglelefteq  G$  and let $\pi$ be a set of 
primes. Then $G\in E_\pi$ if
and only if $G/A\in E_\pi$ and $A$ has a  $\pi$-Hall subgroup $H$
such that $H^A=H^G$.
\end{Cor}

\section{Preliminary results}

The notation of the paper is standard. As we say in Introduction, 
$\pi$  always stands for a set of primes. Given $G$, the set of all $\pi$-Hall
subgroups of $G$ is denoted by  $\operatorname{Hall}_\pi(G)$. The notation
$H\prn  G$ means that $H$ is a pronormal subgroup of~$G$.

\begin{Lemma} \label{base} {\em \cite[Ch. IV, (5.11)]{Suz}}
Let $A$ be a normal subgroup of  $G$. If $H$ is a $\pi$-Hall subgroup of
$G$, then $H \cap A$ is a  $\pi$-Hall subgroup of  $A$, while  $HA/A$ is a
$\pi$-Hall subgroup of~$G/A$.
\end{Lemma}

Recall that a   group is called {\it $\pi$-separable}, if it has a normal
series with factors either  $\pi$- or $\pi'$-groups.

\begin{Lemma} \label{base1} {\em \cite[Ch. V, Theorem
3.7]{Suz}}
Every $\pi$-separable group satisfies~$C_\pi$.
\end{Lemma}

\begin{Lemma}\label{HallExist} {\em \cite[Lemma~2.1(e)]{RevVdoArxive}}
Let $A$ be a normal subgroup of $G$ such that  $G/A$ is a
$\pi$-group, $U$ a $\pi$-Hall subgroup of  $A$. Then a
$\pi$-Hall subgroup  $H$ of $G$ with   $H\cap A=U$ exists if and only if $U^G=U^A$.
\end{Lemma}

\begin{Lemma} \label{Hall_Pron_Simple_Lemma} {\em \cite[Theorem~1]{VR3}} The 
Hall subgroups in simple groups are pronormal.
\end{Lemma}

\begin{Lemma} \label{FrattiniArgument} {\em \cite[Theorem~1]{VR5}}
Let $G\in E_\pi$ and $A\trianglelefteq  G$. Then there exists 
$H\in\operatorname{Hall}_\pi(A)$ such that
$G=AN_G(H)$. Moreover   $N_G(H)\in E_\pi$ and
$\operatorname{Hall}_\pi(N_G(H))\subseteq\operatorname{Hall}_\pi(G)$.
\end{Lemma}

\begin{Lemma} \label{ReplWithConj}
Let $H$ be a subgroup of $G$. Considering $g\in G$, $y\in \langle H,
H^g\rangle$, assume that subgroups  $H^y$ and $H^g$ are conjugate in
$\langle H^y,
H^g\rangle$. Then
$H$ and $H^g$ are conjugate in~$\langle H, H^g\rangle$.
\end{Lemma}

\begin{proof} Let $z\in \langle H^y, H^g\rangle$ and
$H^{yz}=H^g$. Then  $z\in \langle H, H^g\rangle$, since $\langle H^y,
H^g\rangle\le\langle H, H^g\rangle $. So $x=yz\in\langle H, H^g\rangle$ and
$H^x=H^g$.
\end{proof}

The statements of the following two lemmas are evident.

\begin{Lemma} \label{Quot}
Let $\overline{\phantom{g}}:G\rightarrow G_1$ be a group homomorphism and 
$H\leq G$.
Then $H\prn  G$ implies $\overline{H}\prn \overline{G}$.\end{Lemma}

\begin{Lemma} \label{Over} Let $G$ be a group. Then  $H\prn  G$ implies
$H\prn K$, for every subgroup  $K$ of $G$ such that  $H\leq K$.
\end{Lemma}

\begin{Lemma} \label{CentralProd} {\em \cite[Lemma 7]{VR3}}
Let $G$ be a finite group and let  $G_1,\dots,G_n$ be normal subgroups of $G$ 
such
that  $[G_i,G_j]=1$ for
$i\ne j$ and $G=G_1\cdot\ldots\cdot G_n$. Assume that for every  $i=1,\dots, n$
a pronormal subgroup  $H_i$ of
$G_i$ is chosen, and let $H=\langle H_1,\dots, H_n\rangle$.
Then $H\prn  G$.
\end{Lemma}

\begin{Lemma} \label{HallSubgrIn HomImage} {\em \cite[Corollary 9]{ExCrit}}
Let $G\in E_\pi$ and $A\trianglelefteq  G$.
Then for every  $K/A\in\operatorname{Hall}_\pi(G/A)$
there exists $H\in\operatorname{Hall}_\pi(G)$ such that $K=HA$.
\end{Lemma}

\begin{Lemma} \label{ProCrit} Let $H\leq G$ and $A\trianglelefteq  G$. The 
following   are equivalent:
\begin{itemize}
\item[$(1)$] $H\prn  G$.
\item[$(2)$] $HA\prn  G$ and $H\prn  N_G(HA)$.
\end{itemize}
\end{Lemma}

\begin{proof} Assume (1). Then
$HA\prn  G$ by Lemma~\ref{Quot}  and $H\prn  N_G(HA)$ by
Lemma~\ref{Over}. Conversely, assume (2). Take  $g\in G$. We need
to show that there exists $x\in\langle H,H^g\rangle$ such that $H^x=H^g$.
Since
$HA/A\prn  G/A$, there exists $y\in \langle H,H^g\rangle$ with
$H^yA=H^gA$. In accordance with Lemma~\ref{ReplWithConj} it is possible to 
replace  $H$ by
$H^y$, and to assume that  $HA=H^yA=H^gA$, i.~e. $g\in N_G(HA)$.
Since $H\prn  N_G(HA)$, the existence of desired  $x$ is evident.
\end{proof}

\begin{Lemma} \label{PiSep} Let  $A$ be a  $\pi$-separable normal subgroup of $G$ and
${H\in\operatorname{Hall}_\pi(A)}$.
Then $H\prn  G$.
\end{Lemma}

\begin{proof} The lemma follows since the 
subgroup
$\langle H,  H^g\rangle\leq A$ is $\pi$-separable for every  $g\in G$, while if 
$H$ and
$H^g$ are its $\pi$-Hall subgroups, then they are conjugate in   $\langle H, 
H^g\rangle$
by Lemma~\ref{base1}.
\end{proof}

\begin{Lemma} \label{SuffPro}
Let $B$ be a normal subgroup of a finite group  $G$. Then for every normal
subgroup   $A$ of $G$
including $B$, and for every
$H\in\operatorname{Hall}_\pi(A)$ the conditions
\begin{itemize}
\item[$(1)$] $HB/B\prn G/B$;
\item[$(2)$] $(H\cap B)\prn B$;
\item[$(3)$] $(H\cap B)^G=(H\cap B)^B$
\end{itemize}
imply $H\prn  G$..
\end{Lemma}

\begin{proof} Since $HB/B\prn  G/B$,
we have  $HB\prn  G$. By Lemma~\ref{ProCrit}
it is enough to show that $H\prn  N_G(HB)$. So, without loss of
generality, we may assue that $G=N_G(HB)$, i.~e. $HB\trianglelefteq  G$, and
$A=N_A(HB)$. Notice that in this case  $A/B$ is
$\pi$-separable, since the factors $A/HB$ and $HB/B$ of the normal series
$A\trianglerighteq HB\trianglerighteq B$ are   $\pi'$- and $\pi$-groups
respectively.

Take $g\in G$ arbitrary . Since
$(H\cap B)^G=(H\cap B)^B$,  there exists $b\in B$ such that $H^g\cap B=H^b\cap
B$. Since
$(H\cap B)\prn  B$, there exists
$$
y\in \langle H\cap B, H^b\cap B\rangle
=\langle H\cap B, H^g\cap B\rangle\leq \langle H,H^g\rangle
$$
such that  $H^y\cap B=H^b\cap B.$
By Lemma~ \ref{ReplWithConj}   the conjugacy of $H$ and $H^b$ in
$\langle H^,H^b\rangle$ follows from the conjugacy of  $H^y$ and $H^b$ in
$\langle
H^y,H^b\rangle$.
Thus we can replace  $H$ by $H^y$ and suppose that
$$
(H\cap B)^g=(H\cap B)^b=(H\cap B)^y=H\cap B,
$$
i.~e. $g\in N_G(H\cap B)$. It is clear also that  $H\leq
N_A(H\cap B)$.

Since  $(H\cap B)^G=(H\cap B)^B$, applying the Frattini argument we obtain
$G=BN_G(H\cap B)$ and $A=BN_A(H\cap B)$. Therefore,
$$
N_A(H\cap B)/N_B(H\cap B)\cong BN_A(H\cap B)/B=A/B
$$
is $\pi$-separable.
The group $N_B(H\cap B)$ is also  $\pi$-separable, since it has the
normal
$\pi$-Hall subgroup $H\cap B$. Hence,  $N_A(H\cap B)$ is
$\pi$-separable as well. Then from
$H\in\operatorname{Hall}_\pi(N_A(H\cap B))$ by Lemma~\ref{PiSep} applied to $N_G(H\cap
B)$ and its normal
$\pi$-separable subgroup  $N_A(H\cap B)$ we have
$H\prn  N_G(H\cap B)$.
Since
$g\in N_G(H\cap B)$, for some  $x\in \langle H,
H^g\rangle$ the equality $H^x=H^g$ holds. Thus, $H\prn  G$.
\end{proof}

\section{Proof of the main results}

{\it Proof of Theorem} \ref{MainTheorem}. Let $G\in E_\pi$ and  $A\trianglelefteq G$. We need to show that $A$ has a
$\pi$-Hall subgroup  $H$ such that
$H\prn  G$. We proceed by induction on~$|G|$.

If $|G|=1$, we have nothing to prove.

Let $|G|>1$. Choose a minimal normal subgroup  $B$ of $G$ lying in $A$
(note that the inequality $B\not=A$ is not assumed here). Since by Lemma~\ref{base}
$G/B\in E_\pi$, the factor group $A/B$ has a $\pi$-Hall subgroup
$K/B$ such that $K/B\prn  G/B$ by induction. By Lemma~\ref{FrattiniArgument} it 
follows that $B$
possesses a
$\pi$-Hall subgroup $V$ such that $G=BN_G(V)$ or, equivalently,
$V^G=V^B$. This means, in particular, that
$V^K=V^B$ and, by Lemma~6,
there exists $H\in\operatorname{Hall}_\pi(K)$ such that $V=H\cap B$.
Notice that
$|A:H|=|A:K||K:H|$ is a $\pi'$-number, and
so $H\in\operatorname{Hall}_\pi(A)$. Let us show that   $H\prn  G$,
so proving the theorem. We use Lemma~\ref{SuffPro}. By the choice of $K$ we have
$HB/B=K/B \prn  G/B$, which is equivalent to $HB=K\prn
G$, and so (1) of Lemma~\ref{SuffPro} holds. Now  $B$ is a direct product of
simple groups
$
B=S_1\times\dots\times S_n,
$
since $B$ is a minimal normal subgroup of  $G$, and
$
V=\langle V\cap S_i\mid i=1,\dots, n\rangle.
$

Since by Lemma~\ref{base},
$V\cap S_i\in\operatorname{Hall}_\pi(S_i)$
for every $i=1,\dots,n$, and by Lemma ~\ref{Hall_Pron_Simple_Lemma},
$(V\cap S_i)\prn  S_i$, applying Lemma~\ref{CentralProd} we obtain
$
H\cap B=V\prn  B,
$
and so (2) of Lemma~\ref{SuffPro} holds. Finally,
$H\cap B=V$ and by of the choice of $V$ in $B$
we have
$(H\cap B)^G=(H\cap B)^B$. So  (3) of Lemma~\ref{SuffPro} holds.
Thus  $H\prn  G$ by Lemma~\ref{SuffPro}.\qed\medskip

\noindent {\it Proof of Corollary}~\ref{cor}. The   corollary is immediate from 
Theorem~\ref{MainTheorem} for the case~${G=A}$.\qed\medskip

\noindent{\it Proof of Corollary}~\ref{cor1}. Let $A\trianglelefteq  G$. Assume 
that
$G/A\in E_\pi$ and $A$ has a  $\pi$-Hall subgroup  $H$ such that $H^A=H^G$.
Show that $G\in E_\pi$. Choose $X/A\in \operatorname{Hall}_\pi(G/A)$.
Since $A\leq X\leq G$, we have
$$
H^A\subseteq H^X\subseteq H^G;
$$
whence $H^A\subseteq H^X$. Taking into account Lemma~\ref{HallExist} and the fact that
$X/A$ is a  $\pi$-group, we obtain $X\in E_\pi$. Since $|G:X|=|G/A:X/A|$ is a
$\pi'$-number, we have
$$
\varnothing\ne\operatorname{Hall}_\pi(X)
\subseteq\operatorname{Hall}_\pi(G)\ \text{ and }\ G\in E_\pi.
$$

Conersely, let  $G\in E_\pi$. Then by Lemma~\ref{base} $G/A\in E_\pi$. By 
Theorem~\ref{MainTheorem},
there exists a subgroup  $H\in \operatorname{Hall}_\pi(A)$ such that
$H\prn  G$. In particular, for every $g\in G$ there exists
$a\in\langle H, H^g\rangle\leq A$ such that $H^g= H^a$. So  $H^A=H^G$.\qed

\section{Remarks}

\newtheorem{Rem}{{\bfseries Remark}}

\begin{Rem} Theorem~\ref{MainTheorem} guarantees the existence of a pronormal
$\pi$-Hall subgroup in every normal subgroup $A$ of $G\in E_\pi$. The condition
$G\in E_\pi$ cannot be replaced by $A\in
E_\pi$ which is weaker in view of Lemma~\ref{base}. Indeed, let  $\pi=\{2,3\}$ and
$A=\operatorname{GL}_3(2)=\operatorname{SL}_3(2)$.
Then (см.~\cite[Theorem 1.2]{RevHallp}) $A$ has exactly the two classes of 
conjugate  $\pi$-Hall
subgroups with representatives
$$
H_1=\left(
\begin{array}{c@{}c}
\fbox{
$\begin{array}{c}
\!\!\!
\GL_2(2)
\!\!\! \\
\end{array}$
}
& *\\
0 &\fbox{1}
\end{array}
\right){\text{\,\,\, and\,\,\, }}
H_2=\left(
\begin{array}{c@{}c}
\fbox{1}& *\\
 0&\fbox{$\begin{array}{c}
\!
\GL_2(2)
\!
\\
\end{array}$}
\end{array}
\right)
$$
respectively. The class $H_1^A$ consists of line stabilizers in the natural
presentation of
$G$, while $H_2^A$ consists of plane stabilizers. The map  $\iota : x\in
A\mapsto (x^t)^{-1}$ is an automorphism of order  2 of $A$ (here  $x^t$ is the
transpose of~$x$). This automorphism interchanges  $H_1^A$ and $H_2^A$.
Consider the natural split extension
$G=A:\langle\iota\rangle$. Subgroups $H_1$ and $H_2$ are conjugate in $G$. At the
same time $H_1$ and $H_2$ are not conjugate in $A$, cotaining both $H_1$ and $H_2$,
and so they are not pronormal in $G$. We remain to notice that
$G\notin E_\pi$ in this example  by Lemma~\ref{HallExist}.
\end{Rem}

\begin{Rem}
In \cite{VR2,VR4} the definition of   strongly pronormal subgroup is introduced.
Recall that a subgroup  $H$ of $G$ is called {\it strongly pronormal}, if for
every $g\in G$ and $K\leq H$ there exists  $x\in\langle H, K^g\rangle$
such that $K^{gx}\leq H$. In~\cite{VR4} the conjecture that every pronormal Hall subgroup
should be strongly pronormal is formulated. In the light of Theorem~\ref{MainTheorem} and its
corolary it is natural to formulate the weaker conjecture: Does every 
$E_\pi$-group
have a strongly pronormal $\pi$-Hall subgroup?
\end{Rem}

\begin{Rem}
Lemma \ref{SuffPro} plays an important role in the proof of Theorem~\ref{MainTheorem}, and from the lemma,
the following test for pronormality of Hall subgroups follows in particular:
Let $A\trianglelefteq G$ for a group~$G$. A Hall subgroup  $H$ of $G$ is pronormal if

\begin{itemize}
\item[$(1)$] $HA/A\prn G/A$;
\item[$(2)$] $(H\cap A)\prn A$;
\item[$(3)$] $(H\cap A)^G=(H\cap A)^A$.
\end{itemize}

The authors do not know, whether the converse is true; more precisely, whether 
the condition
$H\prn  G$, where  $H$ is a Hall subgroup of $G$, implies 
(2) and (3)? Statement  (1) follows from the pronormality of $H$ by 
Lemma~\ref{Quot}.
\end{Rem}

\end{document}